\magnification=1200

\input amstex

\documentstyle{amsppt}


\hsize=165truemm

\vsize=227truemm


\def\p#1{{{\Bbb P}^{#1}_{k}}}

\def\Hilb{{{\Cal H}\kern -0.25ex{\italic ilb\/}}}

\def\Scand{{{\Cal S}\kern -0.25ex{\italic cand\/}}}

\def\Hom{{{\Cal H}\kern -0.25ex{\italic om\/}}}

\def\Ext{{{\Cal E}\kern -0.25ex{\italic xt\/}}}

\def\Sim{{{\Cal S}\kern -0.25ex{\italic ym\/}}}

\def\Ker{{{\Cal K}\kern -0.25ex{\italic er\/}}}

\def\emdim{\operatorname{emdim}}

\def\char{\operatorname{char}}

\def\lenght{\operatorname{lenght}}

\def\ext{\operatorname{Ext}}

\def\tor{\operatorname{Tor}}

\def\lev{\operatorname{lev}}

\def\mult{\operatorname{mult}}

\def\Soc{\operatorname{Soc}}

\def\M{\operatorname{\frak M}}

\def\N{\operatorname{\frak N}}

\def\ga#1{{\accent"12 #1}}


\topmatter

\title
The Poincar\'e series of a local Gorenstein ring of multiplicity up to 10 is rational
\endtitle

\rightheadtext{rationality of the Poincar\'e series}

\author
Gianfranco Casnati, Roberto Notari
\endauthor

\address
Gianfranco Casnati, Dipartimento di Matematica, Politecnico di Torino,
c.so Duca degli Abruzzi 24, 10129 Torino, Italy
\endaddress

\email
casnati\@calvino.polito.it
\endemail

\address
Roberto Notari, Dipartimento di Matematica \lq\lq Francesco Brioschi\rq\rq, Politecnico di Milano,
via Bonardi 9, 20133 Milano, Italy
\endaddress

\email
roberto.notari\@polimi.it
\endemail



\abstract
Let $R$ be a local, Gorenstein ring with  algebraically closed residue field $k$ of characteristic $0$ and let $P_R(z):=\sum_{p=0}^{\infty}\dim_k(\tor_p^R(k,k))z^p$ be its Poincar\'e series. We compute $P_R$ when $R$ belongs to a particular class defined in the introduction, proving its rationality. As a by--product we prove the rationality of $P_R$ for all local, Gorenstein rings of multiplicity at most $10$.
\endabstract

\endtopmatter

\document

\head
1. Introduction and notation
\endhead
Let $R$ be a Noetherian, local, ring with maximal ideal ${\frak N}$ and residue field $k:=R/{\frak N}$.

The Poincar\'e (or Betti) series of $R$, namely,
$$
P_R(z):=\sum_{p=0}^{\infty}\dim_k(\tor_p^R(k,k))z^p,
$$
has a particular interest in commutative algebra. It is the generating function of the $\tor$--algebra which is a commutative Hopf algebra. Its dual Hopf algebra is isomorphic to the Yoneda $\ext$--algebra of $R$ (see [G--L] for an account on these and other general properties of $P_A$). Moreover the relationships between $P_R$ and some particular quotients of $R$ have been deeply inspected (see e.g. [Ta], [A--L], [G--L], [RR]).

Since J\.P\. Serre conjectured the rationality  of $P_R$ (see e.g. [Se]) and proved it when $R$ is regular, many authors spent their efforts to prove the conjecture. However J. Anick gave an example in [An] of a local ring $R$ having trascendental $P_R$ (see also [Bo] for an example of local, Gorenstein ring with the same property). 

Without any intention of completeness we recall that the conjecture is, in any case, true for several classes of rings $R$: complete intersection ring (see [Ta]), rings $R$ such that $\emdim(R)-2\le{\roman{depth}}(R)$ (see [Sc]), Gorenstein rings with $\emdim(R)-4\le{\roman{depth}}(R)$ (see [Wi] and [J--K--M]), Koszul algebras (see [Mo]).

Let $J\subseteq R$ be a parametric ideal which is a reduction  of $\N$ (see [Ma], Section 14). We can thus consider the local, Artinian ring $A:=R/J$ and its maximal ideal $\M:=\N/J$.
With the above notations, following J\. Sally we say that a ring $R$ is stretched if $\dim_k(\M^2/\M^3)=1$. In [Sa2] the rationality of $P_R$ is proved when $R$ is a stretched Cohen--Macaulay ring. In [E--V2] such a notion has been generalized introducing almost stretched rings, i.e. rings such that $\dim_k(\M^2/\M^3)=2$. Imitating the proof in [Sa2], the authors prove the rationality of $P_R$ also for almost stretched Gorenstein rings.

In this short note, following the argument of [Sa2] (and of [E--V2]), we deal with the next case, namely the case of Gorenstein rings $R$ satisfying, with the above notations, $\dim_k({\M}^2/{\M}^3)=3$ and ${\M}^4=0$. In particular we prove the following

\proclaim{Theorem A}
Let $R$ be a local, Gorenstein, ring satisfying, with the above notations, $\dim_k({\M}^2/{\M}^3)=3$ and ${\M}^4=0$. Assume that the residue field is algebraically closed of characteristic different from $2$, $3$. Then $P_R$ is rational.
\qed
\endproclaim

As a by--product, using all the above mentioned known results about the rationality of $P_R$ and taking into account that passing to $R/J$ does not change the multiplicity (see [Ma], Theorems 14.13 and 17.11), we are finally able to prove also the following

\proclaim{Theorem B}
Let $R$ be a local, Gorenstein, ring with multiplicity at most $10$.  Assume that the residue field is algebraically closed of characteristic $0$. Then $P_R$ is rational.
\qed
\endproclaim

\subhead
Notation
\endsubhead
For the following definitions and results we refer to [Ma].

A ring is a Noetherian ring. Let $R$ be any local ring with maximal ideal ${\frak N}$ and residue field $k:=R/{\frak N}$. We will denote by $\char(k)$ the characteristic of $k$. 

The Samuel function of $R$, $\chi_R(t):=\dim_k(R/{\frak N}^{t+1})$, coincides for large $t$ with a polynomial $p_R(t)$ of degree $\dim(R)$. The multiplicity of $R$ is 
$$
\mult(R):=\lim_{t\to+\infty}{{\dim(R)!}\over {t^{\dim(R)}}}\chi_R(t).
$$
The embedding dimension of $R$ is $\emdim(R):=\chi_R(1)-\chi_R(0)=\dim_k({\frak N}/{\frak N}^2)$.

For a ring $R$ we denote by ${\roman{depth}}(R)$ the maximum lenght of a regular sequence in $\N$.
Recall that $R$ is Cohen--Macaulay if
$\dim(R)={\roman{depth}}(R)$. If, in addition, the injective dimension of $R$ is finite, then $R$ is called Gorenstein. 

Let $A$ be an Artinian ring with maximal ideal $\M$. If $\M^e\ne0$ and $\M^{e+1}=0$ we define the
level of $A$ as $e$ and denote it by $\lev(A)$ (some other authors prefer to use {\sl maximum socle degree}\/ instead of level). In this case the Samuel function coincides with $\lenght(A)$ in the range $t=e,\dots,+\infty$, hence $\mult(A)=\lenght(A)$.

The Hilbert function of $A$, $H_A(t):=\dim_k(\M^t/{\M}^{t+1})$, is the first difference function of $\chi_A$. It is non--zero only in the range $t=0,\dots,e$. We will thus write simply $H_A=(H_A(0),\dots, H_A(e))$. Moreover $\lenght(A)=\sum_{t=0}^{e}H_A(t)$.

The ring $A$ is Cohen--Macaulay and it is Gorenstein if and only if its socle $\Soc(A):=0\colon_A\M$ has dimension $1$ over $k$. In particular $H_A(0)=H_A(e)=1$.

\head
2. Reduction to the Artinian case
\endhead
Let $R$ be a local ring with maximal ideal $\frak N$ and let $r\in {\frak N}$ be a regular element. We recall that in [Ta], it is proved
$$
P_R(z)=\cases
(1+z)P_{R/(r)}(z),&\quad r\in{\frak N}\setminus{\frak N}^2,\\
(1-z^2)P_{R/(r)}(z),&\quad r\in{\frak N}^2.
\endcases
$$
In particular assume now that $R$ is Cohen--Macaulay of positive Krull dimension $h:=\dim(R)$. Then there is a parameter ideal $J\subseteq R$ which is a reduction of $\N$ and such that $A:= R/J$ is a local, Artinian, Cohen--Macaulay ring with maximal ideal $\M:={\frak N}/J$ and $\dim_{k}(A)=\mult(A)=\mult(R)$ (see [Ma], Theorems 14.14, 14.13 and 17.11).

Since $A$ is Artinian, it follows that it is complete, thus $A=S/I$, where $S$ is a regular local ring and $I\subseteq S$ (see [Ma], Theorem 29.4). We can assume that $\emdim(A)=\emdim(S)=\dim(S)$ (i.e. $I$ is contained in the square of the maximal ideal of $S$).
Moreover it is well--known that, if $R$ is Gorenstein, then the same is true for $A$. 

Then, in order to prove the rationality of the Poincar\'e series of the rings $R$ we are considering, it is necessary and sufficient to prove the same property for local, Artinian, Gorenstein rings $A$ with $\dim_k(A)\le10$.

In this case we have the following possibilities. Either $A$ is a complete intersection or $\emdim(A)\le3$ and $A$ is not a complete intersection or $\emdim(A)=4$ or, finally, $H_A$ is one of the following:
$$
(1,n,1,\dots,1),\qquad(1,n,2,1,\dots,1),\qquad (1,5,3,1),
$$
where $n\ge5$. We examine all the above cases but the last one which will be the object of the next section. To this purpose we recall some classical results in the following two theorems.

\proclaim{Theorem 2.1}
Let $\char(k)\ne2$ and let $A$ be a local, Artinian, Gorenstein ring with $\emdim(A)\le4$. Then $P_A$ is rational. 
\endproclaim
\demo{Proof}
If $A$ is a complete intersection see e.g. [Se] and [Ta]. If $A$ is Gorenstein (but not a complete intersection) then the statement has been proved in the remaining cases in [Wi], Theorem 9, when $n=3$, and in [J--K--M], Theorem A, when $n=4$.
\qed
\enddemo

If $H_A=(1,n,1,\dots,1)$, then $A$ is called stretched (see [Sa1]). In [Sa2] the following theorem is proved.

\proclaim{Theorem 2.2}
Let $\char(k)\ne2$ and let $A$ be a local, Artinian, Gorenstein stretched ring. Then $P_A$ is rational.
\qed
\endproclaim

Finally, if $H_A=(1,n,2,1,\dots,1)$, then $A$ is called almost stretched (se [E--V1]). In [E--V2] it is proved the following

\proclaim{Theorem 2.3}
Let $\char(k)=0$ and let $A$ be a local, Artinian, Gorenstein almost stretched ring. Then $P_A$ is rational.
\qed
\endproclaim

Thus, in this section, we have checked Theorem $B$ for each local, Gorenstein, ring $R$ of multiplicity $d\le 10$, except when the quotient $A:=R/J$ of $R$ satisfies $H_A=(1,5,3,1)$.

\remark{Remark 2.4}
In all the above cases the rationality of $P_A$, hence of $P_R$, is proved by computing it explicitely. More precisely, setting $n:=\emdim(A)$, we have that $P_A$ is
$$
{{1}\over{(1-z)^n}}
$$
if $A$ is a complete intersection,
$$
{{(1+z)^3}\over{1-\varepsilon z^2-\varepsilon z^3+z^5}},
$$
where $\varepsilon$ denotes the minimal number of generators of $I$, if $n=3$ and $A$ is Gorenstein (but not a complete intersection),
$$
{{(1+z)^4}\over{f_A(z)}},
$$
if $n=4$ and $A$ is Gorenstein (but not a complete intersection), where $f_A(t)$ is
one of the following according the structure of the $\tor^S\big(A,k)$
$$
\gather
1-\varepsilon z^2-(2\varepsilon-2) z^3-\varepsilon z^4+z^6,\\
1-\varepsilon z^2-(2\varepsilon-5) z^3-(\varepsilon-6) z^4+2z^5-z^6-z^7,\\
1-\varepsilon z^2-(2\varepsilon-2-p) z^3-(\varepsilon-1-2p) z^4+(p+1)z^5-z^7
\endgather
$$
for some $p=1,\dots,\varepsilon$ and $\varepsilon$ as for $n=3$, and finally
$$
{{1}\over{1-nz+z^2}}
$$
if $A$ is either a stretched or almost stretched Gorenstein ring with $n\ge2$. 
\endremark
\medbreak

\head
3. A class of local, Artinian, Gorenstein rings
\endhead
As explained in the previous section we are reduced to study the case of rings $R$ such that the quotient $A:=R/J$ satisfies $H_A=(1,5,3,1)$. So it is necessary to inspect such rings more carefully. In the present section we will deal with local, Artinian, Gorenstein rings $A$ satisfying the more general condition $H_A=(1,n,3,1)$: if $n=2$, then $A$, being Gorenstein, is a complete intersection, thus $\dim_k(\M^3/\M^4)\ge2$ necessarily. Hence the Gorenstein condition forces $n\ge3$.

When $k$ is algebraically closed and $\char(k)\ne2,3$, in Section 4 of [C--N] we gave a complete classification of such rings $A$ assuming the restrictive extra hypothesis that $A$ is also a $k$--algebra. We will recall the method and we will show that such an hypothesis is actually unnecessary.

Let $\M=(a_1,a_2,a_3,\dots,a_n)$. If $a^2\in\M^3$ for each $a\in\M$, then $2a_1a_2=(a_1+a_2)^2-a_1^2-a_2^2\in\M^3$. Thus we can assume $a_1^2$ is an element of a minimal set of generators of $\M^2$. Repeating the above argument for $A/(a_1^2)$ we can finally assume that $\M^2=(a_1^2,a_2^2,a_3^2)$. 

Let $L:=(a_1,a_2,a_3)$ and let $V\subseteq \M/\M^2$ be the corresponding subspace. Thus we have three relations in $I$ of the form
$$
\alpha_1 a_1^2+\alpha_2  a_2^2+\alpha_3a_3^2+2\overline{\alpha}_1 a_2a_3+2\overline{\alpha}_2 a_1a_3+2\overline{\alpha}_3 a_1a_2\in\M^3,\tag 3.1
$$
where $\alpha_i,\overline{\alpha}_j\in
A$, $i,j=1,2,3$, which induce analogous linearly independent relations in $V$. Hence we have a net $\Cal N$ of conics in the projective space ${\Bbb P}(V)$. Let $\Delta$ be the discriminant curve of $\Cal N$ in ${\Bbb P}(V)$. Then $\Delta$ is a plane cubic and the classification of $\Cal N$ depends on the structure of $\Delta$ as explained in [Wa] (notice that the classification described in [Wa] in the complex case, holds as well if the base field $k$ is an algebraically closed field with $\char(k)\ne2,3$).

In what follows we will describe how to modify the results proved in Section 4 of [C--N] for $k$--algebras, in order to extend them to the case of rings. We will examine only one case, namely the case of integral discriminant curve $\Delta$ (see Section 4.2 of [C--N]), the other ones being similar.

In this case we obtain that Relations (3.1) above become
$a_1a_2+a_3^2,a_1a_3,a_2^2-6p a_3^2+q a_1^2\in\M^3$, where $p,q\in A$. We notice that if one of them is in $\M$ then the corresponding monomial is in $\M^3$, thus it can be assumed to be $0$ in what follows. Hence, from now on we will assume that $p,q$ are either invertible or $0$.

We have
$$
\gather
a_1^2a_2=a_1^2a_3=a_1a_2a_3=a_1a_3^2=a_2^2a_3=a_3^3=0,\\
a_2a_3^2=-a_1a_2^2=q a_1^3,\quad
a_2^3=-6pq a_1^3,
\endgather
$$
thus
$\M^2=(a_1^2,a_3^2,a_2a_3)$ and $\M^3=(a_1^3)$. Relation (3.1) thus become
$$
a_1a_2=-a_3^2+\beta_{1,2} a_1^3,\qquad a_1a_3=\beta_{1,3} a_1^3,\qquad a_2^2=\alpha_{2,2}^1 a_1^2+\alpha_{2,2}^3a_3^2+\beta_{2,2} a_1^3,
$$
where $\alpha_{i,j}^h,{\beta}_{i,j}\in A$ are either invertible or $0$, $\alpha_{2,2}^1=-q$, $\alpha_{2,2}^3=6p$.

In general, we have relations of the form
$$
a_ia_j=\alpha_{i,j}^1a_1^2+\alpha_{i,j}^2a_2a_3+\alpha_{i,j}^3a_3^2+\beta_{i,j}a_1^3,\qquad i\ge1,j\ge4,
$$
where the ${\alpha}_{i,j}^h,\beta_{i,j}\in A$ are as above and $\alpha_{i,j}^h=\alpha_{j,i}^h$, $\beta_{i,j}=\beta_{j,i}$.

Via $(a_2,a_3)\mapsto(a_2+\beta_{1,2} a_1^2,a_3+\beta_{1,3} a_1^2)$, we can assume $\beta_{1,2}=\beta_{1,3}=0$.

If $\alpha_{2,2}^1=0$ then $a_2a_3\in\Soc(A)\setminus\M^{3}$, a contradiction. In particular  $\alpha:=\beta_{2,2}a_1-\alpha_{2,2}^1\in A\setminus\M$, whence
$$
a_1a_2=-a_3^2,\qquad
a_1a_3=0,\qquad
a_2^2=-\beta a_1^2+\alpha a_3^2,\tag 3.2
$$
where $\beta$ is invertible and $\alpha:=\alpha_{2,2}^3$.

If $n=3$ we have finished, thus we will assume $n\ge4$ from now on. Via $a_j\mapsto a_j+\alpha_{1,j}^1a_1+\beta_{1,j}a_1^2+\alpha_{2,j}^2a_3+\beta_{2,j}a_3^2+\alpha_{3,j}^2a_2+\beta_{3,j}a_2a_3$, we can assume $\alpha_{1,j}^1=\beta_{1,j}=\alpha_{2,j}^2=\beta_{2,j}=\alpha_{3,j}^2=\beta_{3,j}=0$, $j\ge4$.

Since $a_1a_3=0$, it follows $\alpha_{1,j}^2a_2a_3^2=(a_1a_j)a_3=a_1a_3a_j=(a_3a_j)a_1=\alpha_{3,j}^1a_1^3$, thus $\alpha_{1,j}^2=\alpha_{3,j}^1=0$, $j\ge4$.
Since $a_1a_2a_j=-a_3^2a_j=(a_3a_j)a_3=0$, we have $\alpha_{1,j}^3a_2a_3^2=(a_1a_j)a_2=a_1a_2a_j=(a_2a_j)a_1=\alpha_{2,j}^1a_1^3$, thus $\alpha_{1,j}^3=\alpha_{2,j}^1=0$, $j\ge4$.
Moreover $0=(a_2a_j)a_3=a_2a_3a_j=(a_3a_j)a_2=\alpha_{3,j}^3a_2a_3^2$, thus $\alpha_{3,j}^3=0$, $j\ge4$.
Finally $0=a_2^2a_j=(a_2a_j)a_2=\alpha_{2,j}^3a_2a_3^2$, thus $\alpha_{3,j}^3=0$, $j\ge4$. We conclude that $a_1a_j=a_2a_j=a_3a_j=0$, $j\ge4$.

It follows that $0=(a_1a_i)a_j=(a_ia_j)a_1=\alpha_{i,j}^1a_1^3$, $0=(a_2a_i)a_j=(a_ia_j)a_2=\alpha_{i,j}^3a_2a_3^2$, $0=(a_3a_i)a_j=(a_ia_j)a_3=\alpha_{i,j}^2a_2a_3^2$, thus $a_ia_j=\beta_{i,j}a_1^3$, $i,j\ge4$.

Consider the symmetric matrix $B:=\left(\beta_{i,j}\right)_{i,j\ge4}$ and its class $\overline{B}$ modulo $\M$. If $\det(B)\in\M$, there would exist a non--zero $y:={}^t(y_4,\dots,y_n)\in k^{\oplus n-3}$ such that $\overline{B}y=0$. Let $\alpha_i\in A$ whose residue class in $k$ is $y_i$, $i=4,\dots,n$. By construction the element $a:=\sum_{i=4}^n\alpha_ia_i$ has non--zero class in $\M/\M^2$, thus $a\not\in\M^3$. 

On one hand $a_1a_j=a_2a_j=a_3a_j=0$ (see above). On the other hand
$$
aa_j=\sum_{i=4}^n\alpha_ia_ia_j=(\sum_{i=4}^n\alpha_i\beta_{i,j})a_1^3
$$
if $j\ge4$. Due to the choice of $\alpha_4,\dots,\alpha_n$ we have $\sum_{i=4}^n\alpha_i\beta_{i,j}\in\M$, thus $aa_j=0$, $j\ge4$, hence $a\in\Soc(A)\setminus\M^3$, a contradiction, thus $B$ is invertible. We can thus make a linear change on $a_4,\dots,a_n$ in such a way that
$$
a_ia_j=\delta_{i,j}u_ia_1^3,\qquad i,j\ge4, \tag3.3
$$
where $u_i\in A$ is invertible.

Since $k$ is algebraically closed and $A$ is Artinian, then complete, via a suitable homothety of $(a_1,\dots,a_n)$ we can assume $\beta=u_4=\dots=u_n=1$ due to the following immediate consequence of Hensel's Lemma (see [Ma], Theorem 8.3).

\proclaim{Lemma 3.4}
Let $A$ be an Artinian ring with algebraically closed residue field. If $u\in A$ is invertible, for each positive integer $n\ge1$ there exists $v\in A$ such that $v^n=u$
\qed
\endproclaim

Recall that we are assuming that $A=S/I$, where $S$ is a regular local ring and $I\subseteq S$ is a suitable ideal. Combining Equalities (3.2) and (3.3) above, we have found a minimal set of generators $\{\ x_1,\dots,x_n\ \}$ of the maximal ideal of $S$ such that
$$
(x_1x_2+x_3^2,x_1x_3,x_2^2+x_1^{2}-\alpha x_3^2,x_ix_j,x_h^2-x_1^3)_{{1\le i<j\le n,\ 4\le j}\atop 4\le h\le n}\subseteq I.
$$
where $\alpha\in A$.

\proclaim{Theorem 3.5}
Let $k$, $A$, $S$ and $I$ be as above. Then there is  a minimal set of generators $\{\ x_1,\dots,x_n\ \}$ of the maximal ideal of $S$ such that $I$ is one of the following
$$
\gather
I_1:=(x_1x_2+x_3^2,x_1x_3,x_1^{2}+x_2^2-\alpha x_3^2,x_ix_j,x_h^2-x_1^3)_{{1\le i<j\le n,\ 4\le j}\atop 4\le h\le n},\quad
\alpha\in A,\\
I_{3-p}:=(x_1^2,x_2^2,x_3^2+2px_1x_2,x_ix_j,x_h^2-x_1x_2x_3)_{{1\le i<j\le n,\ 4\le j}\atop 4\le h\le n},\qquad p=0,1,\\
I_4:=(x_2^3-x_1^3,x_3^3-x_1^3,x_ix_j,x_h^2-x_1^3)_{1\le i<j\le n\atop 4\le h\le n},\\
I_5:=(x_1^2,x_1x_2,x_2x_3,x_2^3-x_3^3,x_1x_3^2-x_3^3,x_ix_j,x_h^2-x_3^3)_{{1\le i<j\le n,\ 4\le j}\atop 4\le h\le n},\\
I_6:=(x_1^2,x_1x_2,2x_1x_3+x_2^2,x_3^3,x_2x_3^2,x_ix_j,x_h^2-x_1x_3^2)_{{1\le i<j\le n,\ 4\le j}\atop 4\le h\le n}.
\endgather
$$
\endproclaim
\demo{Proof}
When $\Delta$ is integral we have checked that $I_1\subseteq I$, thus we have an epimorphism $\varphi\colon S/I_1\twoheadrightarrow S/I$, hence $\lenght(S/I_1)\ge\lenght(S/I)$. It is easy to check by direct computation that the Hilbert function of $S/I_1$ is bounded from above by the Hilbert function of $S/I$, hence $\lenght(S/I_1)\le\lenght(S/I)$. We conclude that equality must hold, whence $\varphi$ turns out to be an isomorphism i.e. $I_1=I$.

Thus we have proved above the Theorem when $\Delta$ is integral. When $\Delta$ is not integral the argument is similar: one has to imitate the methods of Section 4.3 of [C--N].
\qed
\enddemo

\head
4. The proof of Theorem A
\endhead
In this section we will prove the rationality of $P_A$ for the rings defined in the previous section, thus completing the proof of Theorem $A$. To this purpose we will imitate again the proof given in [Sa2] (or [E--V2]). We will assume that the residue field of $A$ is algebraically closed with $\char(k)\ne2,3$.

We have to examine the six ideals described in Theorem 3.5. Since the computations are similar, we will consider again only one case, e.g. the last ideal. In this case
$$
I_6:=(x_1^2,x_1x_2,2x_1x_3+x_2^2,x_3^3,x_2x_3^2,x_ix_j,x_h^2-x_1x_3^2)_{{1\le i<j\le n,\ 4\le j}\atop 4\le h\le n},
$$
where $u_4,\dots,u_n$ are invertible.

Notice that in our case the socle of $A$, $\Soc(A)$, is generated by the class  of $x_1x_3^2$. We recall that
$$
P_A(z)={{P_{A/\Soc(A)}(z)}\over{1+z^2P_{A/\Soc(A)}(z)}},\tag4.1
$$
(see [RR] or [A--L]), hence it suffices to prove the rationality of $P_{A/\Soc(A)}(z)$. We  have
$$
A/\Soc(A)\cong A/(x_1x_3^2)\cong S/(x_1^2,x_1x_2,2x_1x_3+x_2^2,x_3^3,x_2x_3^2,x_ix_j,x_h^2,x_1x_3^2)_{{1\le i<j\le n,\ 4\le j}\atop 4\le h\le n}.
$$

Since the classes of $x_4,\dots,x_n$ are trivially in $\Soc(A/(x_1x_3^2))$, it follows from Proposition 3.4.4 of [G--L] that
$$
P_{A/\Soc(A)}(z)={{P_{S_0/J}(z)}\over{1-(n-3)zP_{S_0/J}(z)}}\tag4.2
$$
where $S_0:=S/(x_4,\dots,x_n)$ and $J:=(x_1^2,x_1x_2,2x_1x_3+x_2^2,x_3^3,x_2x_3^2,x_1x_3^2)$, hence it suffices to prove the rationality of $P_{S_0/J}(z)$. Consider now the ideal $H:=(x_1^2,x_1x_2,2x_1x_3+x_2^2,x_3^3,x_2x_3^2)\subseteq S_0$. The quotient $S_0/H$ is a local, Artinian Gorenstein ring since $H$ coincides with the ideal $I_6$ defined in Theorem 3.5 when $n=3$ and its socle is generated by the class of $x_1x_3^2$. Thus again by [RR] or [A--L] we obtain that 
$$
P_{S_0/H}(z)={{P_{S_0/J}(z)}\over{1+z^2P_{S_0/J}(z)}},
$$
hence
$$
P_{S_0/J}(z)={{P_{S_0/H}(z)}\over{1-z^2P_{S_0/H}(z)}}.\tag4.3
$$
Thus it suffices to prove the rationality of $P_{S_0/H}(z)$. 

We finally conclude with the proof of Theorem A in the 

\proclaim{Theorem 4.4}
Let $k$ be algebraically closed with $\char(k)\ne2,3$ and let $A$ be a local, Artinian, Gorenstein, ring with Hillbert function $H_A=(1,n,3,1)$. Then $P_A$ is rational.
\endproclaim
\demo{Proof}
Since $\emdim(S_0/H)=3$, the statement follows by combining Formulas (4.1), (4.2), (4.3) and Theorem 2.2 in the case $A\cong S/I_1$. The same argument gives the rationality of the Poincar\'e series also for $A\cong S/I_t$ with $t=1,\dots,5$.
\qed
\enddemo

As explained in the previous section, an immediate consequence of the above theorem is the proof of Theorem B when $k$ is algebraically closed with $\char(k)=0$.

\remark{Remark 4.5}
Again $P_A$, hence $P_R$, can be explicitely computed via the method described above. Indeed
$$
P_A(z)=\cases{{1}\over{1-nz+3z^2-z^3}}&\text{if $t\le3$},\\
{{(1+z)^3}\over{1-(n-3)z-\varepsilon z^2-\varepsilon z^3+z^5}}&\text{if $t\ge4$},
\endcases
$$
where $\varepsilon={n\choose2}+1$: notice that, also in this case, $\varepsilon$ is the minimal number of generators of $I$.
\endremark
\medbreak

\remark{Remark 4.6}
We list here some comments to our proof and about some possible improvements.

In order to prove Theorem B we have assumed $\char(k)=0$. This is necessary only for proving the rationality of $P_A$ when $A$ is almost stretched. In order to have rings of multiplicity at most $10$ the only interesting cases for the Hilbert function of $A$ are only $(1,5,2,1)$, $(1,5,2,1,1)$ and $(1,6,2,1)$. The other ones are covered by the results in [Ta], [Wi], [J--K--M]. Imitating the method described in Section 3 of [C--N], modified as explained in Section 3 above, one can easily extend the result of rationality of $P_A$ also to the case $\char(k)\ne2,3$ (and $k$ not necessarily algebraically closed).

We now spend a few words about the other condition on $k$, namely algebraic closure. We have not checked all the details but we feel that, via a careful analisys of the classification of nets of conics given in [Wa] which is at the base of the proof of Theorem 3.5, such a condition can be discarded.

If this is true, Theorem A obviously holds under the general hypothesys that the characteristic of $k$ is neither $2$ nor $3$, for each ring $R$ for which there exists a parameter ideal $J\subseteq R$ satisfying $\dim_k(\M^2/\M^3)=3$, where, as usual $\M:=\N/J$ and $\N$ is the maximal ideal of $R$.

Finally we notice that, in order to extend our  result also to all Gorenstein rings of multiplicity $11$, we need to study the Poincar\'e series of local Artinian, Gorenstein rings $A$ satisfying either $H(A)=(1,5,3,1,1)$ or $H(A)=(1,5,4,1)$. This second case seems to be quite difficult since, in our philosophy, it is related to the classification up to projectivities of linear systems of projective dimension $3$ of quadrics in $\p3$.
\endremark
\medbreak

\enddocument